%this is a Plaintex file. The dvi file is created with the instruction
% tex filename
% To tex it successfully you need the macros.tex file

% filename: macros.tex
% definitions used in the TEX model

\def\ifundefined#1#2{\expandafter\ifx\csname#1\endcsname\relax\input #2\fi}
%\ifundefined{Bbb}{amssymb}
\input amssym.def
% load AMS fonts, but only once

% \def\RR{{\rm I \kern-2pt R}}
% \def\CC{{\rm C \kern-6pt C}}
% \def\NN{{\rm I \kern-2pt N}}
% \def\ZZ{{\rm Z \kern-4pt Z}}
% \def\FF{{\rm I \kern-2pt F}}
% \def\QQ{{\rm I \kern-5pt Q}}

\def\NN{{\Bbb N}}
\def\ZZ{{\Bbb Z}}

\def\DD{{\Bbb D}}

\def\Ma{{\cal\char'101}}
\def\Mb{{\cal\char'102}}
\def\Mc{{\cal\char'103}}

\def\Mf{{\cal\char'106}}

\def\Mj{{\cal\char'112}}

\def\Mm{{\cal\char'115}}

\def\Mo{{\cal\char'117}}
\def\Mp{{\cal\char'120}}

\def\Ms{{\cal\char'123}}

\def\pmb#1{\setbox0=\hbox{$#1$}       % generate bold face
     \kern-.025em\copy0\kern-\wd0
     \kern.05em\copy0\kern-\wd0
     \kern-.025em\box0}

%pmb: poor man bold
%to get bold face Phi: $\pmb{\Phi}$

%\def\endproof{$\hfill \diamondsuit$}
\def\endproof{$\hfill \square$}

%\def\ord{\mathop{\rm ord}\nolimits}

%\def\longra#1#2{\,\lower3pt\hbox{${\buildrel\mathop{#1}\over{\lra\atop{#2}}}$}\,}

    % macro with one variable [parameter] \abs x -> |x|
% macro with one variable [parameter] \norm x -> ||x||

\def\Cross{\bigm| \kern-5.5pt \not \ \, }
\def\cross{\mid \kern-5.0pt \not \ \, }             
\def\notto{\hbox{$~\rightarrow~\kern-1.5em\hbox{/}\ \ $}}

\def\Om{\Omega}

\def\sg{\sigma}

\hyphenation{math-ema-ticians}
\hyphenation{pa-ra-meters}
\hyphenation{pa-ra-meter}
\hyphenation{lem-ma}
\hyphenation{lem-mas}
\hyphenation{to-po-logy}
\hyphenation{to-po-logies}
\hyphenation{homo-logy}
\hyphenation{homo-mor-phy}

\def\nSigma{\Sigma \kern-8.3pt \bigm|\,}

\font\teneufm=eufm10
\font\eighteufm=eufm8
\font\fiveeufm=eufm5

\newfam\eufam
\textfont\eufam=\teneufm
\scriptfont\eufam=\eighteufm
\scriptscriptfont\eufam=\fiveeufm

%boxes
\def\boxit#1{\vbox{\hrule\hbox{\vrule\kern2.0pt
       \vbox{\kern2.0pt#1\kern2.0pt}\kern2.0pt\vrule}\hrule}}

%verylongrightarrow{length in points}
\def\vlra#1{\hbox{\kern-1pt
       \hbox{\raise2.38pt\hbox{\vbox{\hrule width#1 height0.26pt}}}
       \kern-4.0pt$\rightarrow$}}

%verylongleftarrow{length in points}
\def\vlla#1{\hbox{$\leftarrow$\kern-1.0pt
       \hbox{\raise2.38pt\hbox{\vbox{\hrule width#1 height0.26pt}}}}}

%verylongdoublearrow{length in points}
\def\vlda#1{\hbox{$\leftarrow$\kern-1.0pt
       \hbox{\raise2.38pt\hbox{\vbox{\hrule width#1 height0.26pt}}}
       \kern-4.0pt$\rightarrow$}}

%longrightarrow{length in points}{formula on top}{formula at bottom}
\def\longra#1#2#3{\,\lower3pt\hbox{${\buildrel\mathop{#2}
\over{{\vlra{#1}}\atop{#3}}}$}\,}

%longleftarrow{length in points}{formula on top}{formula at bottom}
\def\longla#1#2#3{\,\lower3pt\hbox{${\buildrel\mathop{#2}
\over{{\vlla{#1}}\atop{#3}}}$}\,}

%longdoublearrow{length in points}{formula on top}{formula at bottom}
\def\longda#1#2#3{\,\lower3pt\hbox{${\buildrel\mathop{#2}
\over{{\vlda{#1}}\atop{#3}}}$}\,}

\def\overrightharpoonup#1{\vbox{\ialign{##\crcr
	$\rightharpoonup$\crcr\noalign{\kern-1pt\nointerlineskip}
	$\hfil\displaystyle{#1}\hfil$\crcr}}}

\catcode`@=11
\def\@@dalembert#1#2{\setbox0\hbox{$#1\rm I$}
  \vrule height.90\ht0 depth.1\ht0 width.04\ht0
  \rlap{\vrule height.90\ht0 depth-.86\ht0 width.8\ht0}
  \vrule height0\ht0 depth.1\ht0 width.8\ht0
  \vrule height.9\ht0 depth.1\ht0 width.1\ht0 }
\def\dalembert{\mathord{\mkern2mu\mathpalette\@@dalembert{}\mkern2mu}}

\def\@@varcirc#1#2{\mathord{\lower#1ex\hbox{\m@th${#2\mathchar\hex0017 }$}}}
\def\varcirc{\mathchoice
  {\@@varcirc{.91}\displaystyle}{\@@varcirc{.91}\textstyle}%  {\@@varcirc{.57}\scriptstyle}
{\@@varcirc{.45}\scriptscriptstyle}}
\catcode`@=12

\font\tensf=cmss10 \font\sevensf=cmss8 at 7pt
\newfam\sffam
\textfont\sffam=\tensf\scriptfont\sffam=\sevensf

\input amssym.def
\input amssym
\magnification=1200
\font\bigsl=cmsl10 scaled\magstep1   
\font\bigsll=cmsl10 scaled\magstep2
\font\bigslll=cmsl10 scaled\magstep3
\tolerance=500
\overfullrule=0pt
\hsize=6.40 true in
\hoffset=.10 true in
\voffset=0.1 true in
\vsize=8.70 true in
\null
\noindent
{\bf Michigan Math. J. 52 (2004)}
\bigskip\smallskip
\centerline{\bigslll A Purity Theorem for Abelian Schemes} 
\bigskip
\centerline{{\bigsll A}{\sl DRIAN} {\bigsll V}{\sl ASIU}}
\bigskip\bigskip\bigskip\smallskip \centerline{{\bf\bigsll 1. Introduction}} 
\bigskip Let $K$ be the field of fractions of a discrete valuation ring $O$. Let $Y$ be a flat $O$-scheme that is regular, and let
$U$ be an open subscheme of $Y$ whose complement in $Y$ is of
codimension in $Y$ at least 2. We call the pair $(Y,U)$ an extensible pair. Let $q:\Ms\to {\rm Sch}_O$ be a stack over the category ${\rm Sch}_O$ of $O$-schemes endowed with the Zariski topology. Let $\Ms_Z$ be the fibre of $q$ over an $O$-scheme $Z$. Answers to the following Question provide
information on $\Ms$.  
\medskip\noindent
{\it {\bigsll Q}}UESTION {\bigsl 1.1.} Is the pull-back functor $\Ms_Y\to \Ms_U$ surjective on  objects?
\medskip 
 Question 1.1 has a positive answer in any one of the following three cases:
\medskip
{\it (i)} $\Ms$ is the stack of morphisms into the N\`eron model over $O$ of an abelian variety over $K$, and $Y$ is smooth over $O$ (see [N]); 
\smallskip
{\it (ii)} $\Ms$ is the stack of smooth, geometrically connected, projective curves of genus at least 2 (see [M-B]);
\smallskip
{\it (iii)} $\Ms$ is the stack of stable curves of locally constant type, and there is a divisor $DIV$ of $Y$ with normal crossings such that the reduced scheme $Y\setminus U$ is a closed subscheme of $DIV$ (see [dJO]). 
\smallskip
Let $p$ be a prime. If the field $K$ is of characteristic 0, then an example of Raynaud--Gabber--Ogus shows that Question 1.1 does not always have a positive answer if $\Ms$ is the stack of abelian schemes (see [dJO, Sec. 6]). This invalidates [FaC, Chap. IV, Thms. 6.4, $6.4^\prime$, 6.8] and leads to the following problem.
\medskip\noindent
{\it {\bigsll P}}ROBLEM {\bigsl 1.2.} Classify all those $Y$ with the property that, for any extensible pair $(Y,U)$ with $U$ containing $Y_K$, every abelian scheme (resp., every $p$-divisible group) over $U$ extends to an abelian scheme (resp., to a $p$-divisible group) over $Y$.
\medskip
We call such $Y$ a healthy (resp., $p$-healthy) regular scheme (cf. [V, 3.2.1(2),(9)]). The counterexample of [FaC, p. 192] and the classical purity theorem of [G, p. 275] indicate that Problem 1.2 is of interest only if $K$ is of characteristic $0$ (resp., only if $O$ is a faithfully flat $\ZZ_{(p)}$-algebra). We shall therefore assume hereafter that $O$ is of mixed characteristic $(0,p)$. Let $e\in\NN$ be the index of ramification of $O$. If $e\le p-2$, then a result of Faltings states that $Y$ is healthy and $p$-healthy regular, provided it is formally smooth over $O$ (see [Mo, 3.6] and [V, 3.2.2(1) and 3.2.17], a correction to step B of which is implicitly achieved here by Proposition 4.1). If $p\ge 5$, then there are local $O$-schemes that are healthy and $p$-healthy regular but are not formally smooth over some discrete valuation ring (see [V, 3.2.2(5)]). The goal of this paper is to prove the following theorem.
\medskip\noindent
{\it {\bigsll T}}HEOREM {\bigsl 1.3.} {\it If $e=1$, then any regular, formally smooth $O$-scheme is healthy and $p$-healthy regular.}
\medskip
The case $p\ge 3$ is already known, as remarked previously. The case $p=2$ answers a question of Deligne. In Section 2 we present complements on the crystalline contravariant Dieudonn\'e functor. These complements are needed in Section 3 to prove Lemma 3.1, which pertains to extensions of short exact sequences of finite, flat, commutative group schemes. In Section 4 we use Lemma 3.1 and [FaC] to prove Theorem 1.3.
\smallskip 
Milne used an analogue of Question 1.1(i) to define integral canonical models of Shimura varieties (see [Mi, Sec. 2] and [V, 3.2.3, 3.2.6]). Theorem 1.3 implies the uniqueness of such integral canonical models and extends parts of [V] to arbitrary mixed characteristic (see [V, 3.2.3.2, 3.2.4, 3.2.12, etc.]). Also one can use Theorem 1.3 and the integral models of compact, unitary Shimura varieties used in [K] to provide the first concrete examples of N\`eron models (as defined in [BLR, p. 12]) of projective varieties over $K$ whose extensions to $\overline{K}$ are not embeddable into abelian varieties over $\overline{K}$. 
\medskip\noindent
{\it {\bigsll A}CKNOWLEDGEMENTS.} We would like to thank G. Faltings for mentioning to us that his result should hold for $p=2$ as well. I also thank W. McCallum and the referee for several suggestions.  
\bigskip\smallskip
\centerline{{\bf \bigsll\bigsll 2. The Crystalline Dieudonn\'e Functor}}
\bigskip
Let $k$ be a perfect field of characteristic $p>0$. Let $\sg_k$ be the Frobenius automorphism of the Witt ring $W(k)$ of $k$, and let $R$ be a regular, formally smooth $W(k)$-algebra. Let $Y:={\rm Spec}(R)$. Let $\Phi_R$ be a Frobenius lift of the $p$-adic completion $R^\wedge$ of $R$ that is compatible with $\sg_k$. Let $\Om_R^\wedge$ be the $p$-adic completion of the $R$-module of relative differentials of $R$ with respect to $W(k)$, and let $d\Phi_{R/p}$ be the differential of $\Phi_R$ divided by $p$. For $n\in\NN$, the reduction mod $p^n$ of $d\Phi_{R/p}$ is denoted in the same way. If $Z$ is an arbitrary $\ZZ_{(p)}$-scheme, let
$$p-FF(Z)$$ 
be the category of finite, flat, commutative group schemes of $p$-power order over $Z$. 
\smallskip
Let $\Mm\Mf_{[0,1]}^\nabla(Y)$ be the Faltings--Fontaine category defined as follows. Its objects are quintuples 
$$(M,F,\Phi_0,\Phi_1,\nabla),$$
where $M$ is an $R$-module, $F$ is a direct summand of $M$, both $\Phi_0:M\to M$ and $\Phi_1:F\to M$ are $\Phi_R$-linear maps, and $\nabla:M\to M\otimes_R \Om_R^\wedge$ is an integrable, nilpotent mod $p$ connection on $M$, such that the following five axioms hold:
\medskip
{\it 1.} $\Phi_0(m)=p\Phi_1(m)$ for all $m\in F$;
\smallskip
{\it 2.} $M$ is $R$-generated by $\Phi_0(M)+\Phi_1(F)$;
\smallskip
{\it 3.} $\nabla\circ\Phi_0(m)=p(\Phi_0\otimes d\Phi_{R/p})\circ\nabla(m)$ for all $m\in M$;
\smallskip
{\it 4.} $\nabla\circ\Phi_1(m)=(\Phi_0\otimes d\Phi_{R/p})\circ\nabla(m)$ for all $m\in F$; and
\smallskip
{\it 5.} locally in the Zariski topology of $Y$, $M$ is a finite direct sum of $R$-modules of the form $R/p^sR$, where $s\in\NN\cup\{0\}$.
\medskip
A morphism $f:(M,F,\Phi_0,\Phi_1,\nabla)\to (M^\prime,F^\prime,\Phi_0^\prime,\Phi_1^\prime,\nabla^\prime)$ between two such quintuples is an $R$-linear map $f_0:M\to M^\prime$ taking $F$ into $F^\prime$ and such that the following three identities hold: $\Phi_0^\prime\circ f_0=f_0\circ\Phi_0$, $\Phi_1^\prime\circ f_0=f_0\circ\Phi_1$ and $\nabla^\prime\circ f_0=(f_0\otimes_R 1_{\Om_R^\wedge})\circ\nabla$. We refer to $M$ as the underlying $R$-module of $(M,F,\Phi_0,\Phi_1,\nabla)$. Disregarding the connections (and thus axioms 3 and 4), we obtain the category $\Mm\Mf_{[0,1]}(Y)$. Categories like $\Mm\Mf_{[0,1]}(Y)$ and $\Mm\Mf_{[0,1]}^\nabla(Y)$, in the context of arbitrary smooth $W(k)$-schemes, were first introduced in [Fa] as inspired by [F] and [FL], which worked with the category $\Mm\Mf_{[0,1]}({\rm Spec}(W(k)))$. In the sequel we will need the following result of Faltings. 
\medskip\noindent
{\it {\bigsll P}}ROPOSITION {\bigsl 2.1.} {\it We assume that $\Om_R^\wedge$ is a flat $R$-module. Then the category $\Mm\Mf_{[0,1]}^\nabla(Y)$ is abelian and the functor from it into the category of $R$-modules that takes $f$ into $f_0$ is exact.}
\medskip\noindent
{\it Proof.} This follows from [Fa, pp. 31--33]. Strictly speaking, in [Fa] the result is stated only for smooth $W(k)$-algebras, but the inductive arguments work also for regular, formally smooth $W(k)$-algebras. In fact, we can use Artin's approximation theorem to reduce Proposition 2.1 to the result in [Fa] as follows. 
\smallskip
Let $f$ and $f_0$ be as before. We denote also by $\Phi_0$, $\Phi_1$, $\nabla$ and $\Phi_0^\prime$, $\Phi_1^\prime$, $\nabla^\prime$ the different $\Phi_R$-linear maps and connections obtained from them via restrictions or via natural passage to quotients (for $\nabla$ and $\nabla^\prime$ this makes sense because $\Om_R^\wedge$ is a flat $R$-module). We need to show that the three quintuples $({\rm Ker}(f_0),F\cap {\rm Ker}(f_0),\Phi_0,\Phi_1,\nabla)$, $(f_0(M),f_0(F),\Phi_0^\prime,\Phi_1^\prime,\nabla^\prime)$ and $(M^\prime/f_0(M),F^\prime/f_0(F),\Phi_0^\prime,\Phi_1^\prime,\nabla^\prime)$ are objects of $\Mm\Mf_{[0,1]}^\nabla(Y)$ and that $f_0(F)=F^\prime\cap f_0(M)$. Since $\Om_R^\wedge$ is a flat $R$-module, axioms 3 and 4 hold and so from now on we do not mention $\nabla$ and $\nabla^\prime$. Hence we are interested only in the morphism $g:(M,F,\Phi_0,\Phi_1)\to (M^\prime,F^\prime,\Phi_0^\prime,\Phi_1^\prime)$ of $\Mm\Mf_{[0,1]}(Y)$ defined by $f_0$. We can assume that $M$ and $M^\prime$ are annihilated by $p^n$ and that $R$ is local. Using devissage as in [Fa, p. 33, ll. 4--11], it is enough to handle the case $n=1$. So all the $R$-modules involved in the three quintuples listed are in fact $R/pR$-modules. Thus, to check that they are free, we can also assume that $R$ is complete. Based on [Ma, p. 268], there is a $k$-subalgebra $k_1$ of $R/pR$ that is isomorphic to the residue field of $R$. We easily get that $R/pR$ is a $k$-algebra of the form $k_1[[x_1,...,x_d]]$, where $d\in\NN\cup\{0\}$. Because $n=1$, the choice of $\Phi_R$ plays no role in the study of the three quintuples and so we can  also assume that $k_1$ is perfect. 
\smallskip
We choose $R/pR$-bases $\Mb$ and $\Mb^\prime$ of $M$ and $M^\prime$ (respectively) such that certain subsets of them are $R/pR$-bases of $F$ and $F^\prime$. With respect to $\Mb$ and $\Mb^\prime$, the functions $f_0$, $\Phi_0$, $\Phi_1$, $\Phi_0^\prime$, and $\Phi_1^\prime$ involve a finite number of coordinates that are elements of $R/pR$. Let $A_0$ be the  $k_1$-subalgebra of $R/pR$ generated by all these coordinates, and observe that $A_0$ is of finite type. Hence, from [BLR, p. 91] we derive the existence of an $A_0$-algebra $A_1$ that is smooth over $k_1$ and such that the $k_1$-monomorphism $A_0\hookrightarrow R/pR$ factors through $A_1$. Localizing $A_1$, we can assume that $A_1$ is the reduction mod $p$ of a smooth $W(k_1)$-algebra $R_1$. Now fix a Frobenius lift of the $p$-adic completion of $R_1$ that is compatible with $\sg_{k_1}$; hence we can speak about $\Mm\Mf_{[0,1]}(R_1)$. We get that $g$ is the natural tensorization with $R$ of a morphism $g_1$ of $\Mm\Mf_{[0,1]}(R_1)$. Applying [Fa, pp. 31--32] to $g_1$ and tensoring with $R$, we deduce that axioms 1, 2, and 5 hold for the three quintuples and that $f_0(F)=F^\prime\cap f_0(M)$.\endproof 
\medskip\noindent
{\it {\bigsll C}}ONSTRUCTION {\bigsl 2.2.} Let $W_n(k):=W(k)/p^{n}W(k)$. There is a contravariant, $\ZZ_p$-linear functor 
$$\DD:p-FF(Y)\to \Mm\Mf_{[0,1]}^\nabla(Y).$$ 
Similar functors but with $Y$ replaced by ${\rm Spec}(W(k))$ (resp., by a smooth $W(k)$-scheme and with $p>2$) were first considered in [F] (resp. [Fa]). The existence of $\DD$ is a modification of a particular case of [BBM, Chap. 3]. We now include the construction of $\DD$ based in essence on [BBM] and [Fa, 7.1]. We will use Berthelot's crystalline site ${\rm CRIS}(Y_{W_n(k)}/{\rm Spec}(W(k)))$ (see [B, Chap. III, Sec. 4]) and its standard exact sequence $0\to\Mj_{Y_{W_n(k)}/W(k)}\to\Mo_{Y_{W_n(k)}/W(k)}$ (see [BBM, p. 12]). 
\smallskip
Let $G$ be an object of $p-FF(Y)$ that is annihilated by $p^n$.  Let $(\tilde M,\tilde {\Phi}_0,\tilde{V}_0,\tilde{\nabla})$ be the evaluation of the Dieudonn\'e crystal $\DD(G_{Y_{k}})=Ext^1_{Y_k/W(k)}(\underline{G_{Y_k}},\Mo_{Y_k/W(k)})$ (see [BBM, p. 116]) at the thickening naturally attached to the closed embedding $Y_{k}\hookrightarrow Y_{W_n(k)}$. Hence $\tilde M$ is an $R$-module, $\tilde{\Phi}_0$ is a $\Phi_R$-linear endomorphism of $\tilde M$, $\tilde{V}_0:\tilde M\to \tilde M\otimes_{R}\, _{\Phi_R}R$ is a Verschiebung map, and $\tilde\nabla$ is an integrable and nilpotent mod $p$ connection on $\tilde M$. Identifying $\tilde{\Phi}_0$ with an $R$-linear map $\tilde M\otimes_R\, _{\Phi_R}R\to \tilde M$, we have 
$$\tilde{V}_0\circ\tilde{\Phi}_0(x)=px\;\;\;\forall x\in \tilde M\otimes_{R}\, _{\Phi_R}R,\;\;\;{\rm and}\;\;\;\tilde{\Phi}_0\circ\tilde{V}_0(x)=px\;\;\;\forall x\in \tilde M.\leqno (1)$$ 
Let $\tilde F$ be the direct summand of $\tilde M$ that is the Hodge filtration defined by the lift $G_{Y_{W_n(k)}}$ of $G_{Y_k}$. The triple $(\tilde M,\tilde\Phi_0,\tilde V_0,\tilde\nabla)$ is also the evaluation of $\DD(G_{Y_{W_n(k)}})=Ext^1_{Y_{W_n(k)}/W(k)}(\underline{G_{Y_{W_n(k)}}},\Mo_{Y_{W_n(k)}/W(k)})$ at the trivial thickening of $Y_{W_n(k)}$. So $\tilde F$ is the image of the evaluation at this trivial thickening of the functorial  homomorphism $Ext^1_{Y_{W_n(k)}/W(k)}(\underline{G_{Y_{W_n(k)}}},\Mj_{Y_{W_n(k)}/W(k)})\to Ext^1_{Y_{W_n(k)}/W(k)}(\underline{G_{Y_{W_n(k)}}},\Mo_{Y_{W_n(k)}/W(k)})$. 
\smallskip
To define the map $\tilde{\Phi}_1:\tilde F\to\tilde M$ and to check that axioms 1--5 hold for the quintuple $(\tilde M,\tilde F,\tilde{\Phi}_0,\tilde{\Phi}_1,\tilde{\nabla})$, we can work locally in the Zariski topology of $Y$. Hence we can assume that $G$ is a closed subgroup of an abelian scheme $A^\prime$ over $Y$ (cf. Raynaud's theorem of [BBM, 3.1.1]). Let $A:=A^\prime/G$, and let $i_G:A^\prime\twoheadrightarrow A$ be the resulting isogeny. We now define $\tilde{\Phi}_1$ using the cokernel of a morphism $f$ of $\Mm\Mf_{[0,1]}^\nabla(Y)$ associated naturally to $i_G$. 
\smallskip
Let $R(n):=R/p^nR$. Let $M:=H^1_{\rm crys}(A_{R(n)}/R(n))=H^1_{\rm dR}(A_{R(n)}/R(n))$ as in  [BBM, 2.5]. Let $F$ be the direct summand of $M$ that is the reduction mod $p^n$ of the Hodge filtration $F_A$ of 
$$H^1_{\rm crys}(A/R^\wedge):={\rm proj.lim.}_{l\in\NN} H^1_{\rm crys}(A_{R(l)}/R(l))={\rm proj.lim.}_{l\in\NN} H^1_{\rm dR}(A_{R(l)}/R(l)).$$ 
Now let $\Phi_0$ be the reduction mod $p^n$ of the $\Phi_R$-linear endomorphism $\Phi_A$ of $H^1_{\rm crys}(A/R^\wedge)$, and let $\Phi_1$ be the reduction mod $p^n$ of the $\Phi_R$-linear map $F_A\to H^1_{\rm crys}(A/R^\wedge)$ taking $m\in F_A$ into $\Phi_A(m)/p$. Let $\nabla$ be the reduction mod $p^n$ of the Gauss--Manin connection $\nabla_A$ of $A_{R^\wedge}$. That $\Mc:=(M,F,\Phi_0,\Phi_1,\nabla)$ is an object of $\Mm\Mf_{[0,1]}^\nabla(Y)$ is implied by the fact that the quadruple $(H^1_{\rm crys}(A/R^\wedge),F_A,\Phi_A,\nabla_A)$ is the evaluation at the thickening attached naturally to the closed embedding $Y_k\hookrightarrow Y^\wedge:={\rm Spec}(R^\wedge)$ of a filtered $F$-crystal over $R/pR$ in locally free sheaves (see [Ka, Sec. 8]). Similarly, starting from $A^\prime$ we construct $\Mc^\prime=(M^\prime,F^\prime,\Phi_0^\prime,\Phi_1^\prime,\nabla^\prime)$. Let $f:\Mc\to\Mc^\prime$ be the morphism of $\Mm\Mf_{[0,1]}^\nabla(Y)$ associated naturally to $i_G$. 
\smallskip
Let $f_0:M\to M^\prime$ defining $f$. Let 
$$\DD(G)=(\tilde M,\tilde F,\tilde{\Phi}_0,\tilde{\Phi}_1,\tilde{\nabla}):={\rm Coker}(f)$$ 
(cf. Proposition 2.1). Then $\tilde M:=M^\prime/f_0(M)$, $\tilde F:=F^\prime/f_0(F)$, and so forth. That the quadruple $(\tilde M,\tilde F,\tilde{\Phi}_0,\tilde{\nabla})$ is as defined previously follows from [BBM, 3.1.6, 3.2.9, 3.2.10]. 
\smallskip
The association $G\to (\tilde M,\tilde F,\tilde{\Phi}_0,\tilde{\nabla})$ is functorial. In order to check that $\tilde{\Phi}_1$ is well-defined and functorial, we can assume that $R$ is local. To ease the notations we will check directly that $\DD(G)$ is itself well defined and functorial. So let $m:G\to H$ be a morphism of $p-FF(Y)$. If $H$ is a closed subgroup of an abelian scheme $B^\prime$ over $R$, then $\DD(G\times_Y H)$ is computed via the product embedding of $G\times_Y H$ into $A^\prime\times_Y B^\prime$. We thus obtain $\DD(G)\oplus \DD(H)=\DD(G\times_Y H)$. We now define $\DD(m)$. If $m$ is a closed embedding, then the construction of $\DD(m)$ is obvious because $i_G$ factors through the isogeny $i_H:A^\prime\to A^{\prime}/H$. In general, the homomorphism $(1_G,m):G\to G\times_Y H$ is a closed embedding. Hence $\DD(m):\DD(H)\to\DD(G)$ is defined naturally via the epimorphism $\DD(1_G,m):\DD(G)\oplus \DD(H)=\DD(G\times_Y H)\twoheadrightarrow\DD(G)$.
\smallskip
One easily checks that $\DD(G)$ and $\DD(m)$ are well-defined; that is, they depend neither on the chosen embeddings into abelian schemes nor on the choice of a power of $p$ annihilating $G$ and $H$. For instance, let $G$ be a closed subgroup of another abelian scheme $C^\prime$ over $Y$. By embedding $G$ diagonally into $A^\prime\times_Y C^\prime$ and then using the snake lemma in the context of any one of the two projections of $A^\prime\times_Y C^\prime$ onto its factors, we get that $\DD(G)$ defined via $A^\prime\times_Y C^\prime$ is isomorphic to $\DD(G)$ defined via $A^\prime$ or $C^\prime$. This ends the construction of $\DD$.
\medskip\noindent
{\it {\bigsll R}}EMARKS {\bigsl 2.3.} {\it (1)} We have 
$$\tilde{V}_0\circ\tilde{\Phi}_1(x)=x\;\;\;\forall x\in \tilde F\otimes_{R}\, _{\Phi_R}R,\leqno (2)$$
 as this identity holds in the context of $A$ and $A^\prime$. Since $\tilde M$ is $R$-generated by the images of $\tilde{\Phi}_1$ and $\tilde{\Phi}_0$, it follows that $\tilde{V}_0$ is uniquely determined by $\tilde{\Phi}_0$ and $\tilde{\Phi}_1$. We therefore deem it appropriate to denote $(\tilde M,\tilde F,\tilde{\Phi}_0,\tilde{\Phi}_1,\tilde{\nabla})$ by $\DD(G)$. As $\Mc$ and $\Mc^\prime$ depend only on $A_{Y_{W_{n+1}(k)}}$ and $A^{\prime}_{Y_{W_{n+1}(k)}}$ (respectively), $\DD(G)$ also depends only on $G_{Y_{W_{n+1}(k)}}$.
\smallskip
{\it (2)} If $\tilde F$ is neither $\{0\}$ nor $\tilde M$, then $\tilde{V}_0$ has a nontrivial kernel and so $\tilde{\Phi}_1$ is not determined by $\tilde{V}_0$. The advantage we gain by using $\tilde{\Phi}_1$ instead of $\tilde{V}_0$ is that we can exploit axiom 5 and the exactness part of Proposition  2.1 (see the proof of Lemma 3.1).
\smallskip
{\it (3)} Let $Y_1={\rm Spec}(R_1)$ be an affine, regular, formally smooth $W(k)$-scheme. We assume that $R_1^\wedge$ is equipped with a Frobenius lift $\Phi_{R_1}$ compatible with $\sg_k$ and that there is a morphism $l:Y_1\to Y$ whose $p$-adic completion $l^\wedge$ is compatible with the Frobenius lifts. Let $l^*:p-FF(Y)\to p-FF(Y_1)$ and $l^*:\Mm\Mf_{[0,1]}^\nabla(Y)\to\Mm\Mf_{[0,1]}^\nabla(Y_1)$ be the pull-back functors. Hence $l^*(G)=G\times_Y {Y_1}$ and 
$$l^*(M,F,\Phi_0,\Phi_1,\nabla)=(M\otimes_R R_1,F\otimes_R R_1,\Phi_0\otimes\Phi_{R_1},\Phi_1\otimes\Phi_{R_1},\nabla_1),$$ 
where $\nabla_1$ is the natural extension of $\nabla$ to a connection on $M\otimes_R R_1$. These constructions then yield the equality $\DD\circ l^*=l^*\circ\DD$ of contravariant, $\ZZ_p$-linear functors from $p-FF(Y)$ to $\Mm\Mf_{[0,1]}^\nabla(Y_1)$. 
\smallskip
{\it (4)} As in [Fa, 2.3], we see that the category $\Mm\Mf_{[0,1]}^\nabla(Y)$ does not depend (up to isomorphism) on the choice of the Frobenius lift $\Phi_R$ of $R^\wedge$ compatible with $\sg_k$. The arguments of [Fa] apply even for $p=2$ because we are dealing with connections that are nilpotent mod $p$. One can use this to show that remark (3) makes sense even if $Y$ and $Y_1$ are not affine or if no Frobenius lifts are fixed. 
\smallskip
{\it (5)} If $R$ is local, complete, and has residue field $k$, then one can use a theorem of Badra [Ba] on the category $p-FF(Y)$ to obtain directly that $\DD(G)$ is functorial.
\bigskip\smallskip
\centerline{{\bf \bigsll\bigsll 3. A Lemma}}
\bigskip 
In this section we prove the following Lemma. 
\medskip\noindent
{\it {\bigsll L}}EMMA {\bigsl 3.1.} {\it Assume that $e=1$. Let $(Y,U)$ be an extensible pair, with $Y$ a regular and formally smooth $O$-scheme of dimension $2$ and with $U$ containing $Y_K$. Then any short exact sequence $0\to G_{1U}\to G_{2U}\to G_{3U}\to 0$ in the category $p-FF(U)$ extends uniquely to a short exact sequence in the category $p-FF(Y)$.}
\medskip\noindent
{\it Proof.} Let $\Mo_X$ be the sheaf of rings on a scheme $X$. Let $j:U\hookrightarrow Y$ be the open embedding of $U$ in $Y$. For $i\in\{1,2,3\}$, the $\Mo_Y$-module $\Mf_i:=j_*(\Mo_{G_{iU}})$ is locally free (cf. [FaC, Lemma 6.2 of p. 181]). The commutative Hopf algebra structure of the $\Mo_{U}$-module $\Mo_{G_{iU}}$ extends uniquely to a commutative Hopf algebra structure of $\Mf_i$. Hence there exists a unique finite, flat, commutative group scheme $G_i$ over $Y$ extending $G_{iU}$. We have to show that the natural complex
$$0\to G_1\to G_2\to G_3\to 0\leqno (3)$$
is, in fact, a short exact sequence. This is a local statement for the faithfully flat topology of $Y$. We may therefore assume that $Y$ is local and complete and that its residue field $k$ is separable closed and of characteristic $p$; we may also assume that $U$ is the complement in $Y$ of the maximal point $y$ of $Y$. We write $Y={\rm Spec}(R)$. From Cohen's coefficient ring theorem (see [Ma, pp. 211, 268]) we have that $R$ is a $K(k)$-algebra, where $K(k)$ is a Cohen ring of $k$. Since $R/pR$ is regular and formally smooth over $O/pO$ (and thus also over $k$),  we can identify $R=K(k)[[x]]$ as $K(k)$-algebras. Hence, by replacing $R$ with the faithfully, flat $R$-algebra $W(\overline{k})[[x]]$, we can assume that $k=\overline{k}$ and $K(k)=W(k)$ and so can use the notations of Section 2 (e.g. $\Phi_R$, $\Om_R^\wedge$, ...). Since $\Om_R^\wedge=dxR$ is a free $R$-module, we can also appeal to Proposition 2.1.
\smallskip
Let $\Mo$ be the local ring of $Y$, which is a discrete valuation ring that is faithfully flat over $W(k)$. Let $\Mo_1:=W(k_1)$, where $k_1$ is the algebraic closure of the residue field $k((x))$ of $\Mo$. We consider a Teichm\"uller lift $l:{\rm Spec}(\Mo_1)\to {\rm Spec}(R^\wedge)$ that---at the level of special fibres---induces the inclusion $k[[x]]\hookrightarrow k_1$. Hence, $\Mo_1$ has a natural structure of an $\Mo$-algebra. Let 
$$0\to\DD(G_3)\to\DD(G_2)\to\DD(G_1)\to 0\leqno (4)$$
be the complex of $\Mm\Mf_{[0,1]}^\nabla(Y)$ corresponding to (3). Let $M_1$, $M_2$ and $M_3$ be the underlying $R$-modules of $\DD(G_1)$, $\DD(G_2)$, and $\DD(G_3)$, respectively. Let
$$0\to M_3\to M_2\to M_1\to 0\leqno (5)$$
be the complex of $R$-modules defined by (4). Let $N_{1,2}$ be the underlying $R$-module of ${\rm Coker}(\DD(G_2)\to\DD(G_1))$. The key point is that ${\rm Coker}(\DD(G_2)\to\DD(G_1))$ exists in the category $\Mm\Mf^{\nabla}_{[0,1]}(Y)$ and the sequence $M_2\to M_1\to N_{1,2}\to 0$
is exact (cf. Proposition 2.1). We show that $N_{1,2}=\{0\}$. Because $N_{1,2}$ is a direct sum of $R$-modules of the form $R/p^sR=W_s(k)[[x]]$ for $s\in\NN\cup\{0\}$ (cf. axiom 5), to show that $N_{1,2}=\{0\}$ it is enough to show that $N_{1,2}[{1\over x}]=\{0\}$. It is thus enough to show that the complex 
$$0\to M_3\otimes_{\Mo} \Mo_1\to M_2\otimes_{\Mo} \Mo_1\to M_1\otimes_{\Mo} \Mo_1\to 0\leqno (6)$$
 obtained from (5) by tensoring with $\Mo_1$ is a short exact sequence. Note that (6) is the complex obtained by pulling back (3) to ${\rm Spec}(\Mo_1)$, applying $\DD$, and then taking underlying $\Mo_1$-modules (cf. Remark 2.3(4) applied to $l$). But the pull-back of (3) to ${\rm Spec}(\Mo_1)$ is a short exact sequence (since the pull-back of (3) to $U$ is so). Thus (6) is the complex associated via the classical contravariant Dieudonn\'e functor to the short exact sequence $0\to G_{1k_1}\to G_{2k_1}\to G_{3k_1}\to 0$ (cf. [BBM, pp. 179--180]). From the classical Dieudonn\'e theory we threfore have that (6) is a short exact sequence, cf. [F, p. 128 or p. 153]. So $N_{1,2}=\{0\}$. 
\smallskip
Hence the natural $W(k)$-linear map $j_{1,2}:M_2/(x)M_2\to M_1/(x)M_1$ is an epimorphism. But $j_{1,2}$ is the $W(k)$-linear map associated via the classical contravariant Dieudonn\'e functor to the homomorphism $G_{1k}\to G_{2k}$, so this homomorphism is a closed embedding (cf. the classical Dieudonn\'e theory). It follows by Nakayama's lemma that  $G_1$ is a closed subgroup of $G_2$. Both $G_3$ and $G_2/G_1$ are finite, flat, commutative group schemes extending $G_{3U}$ and so we have $G_3=G_2/G_1$. Hence (3) is a short exact sequence. This completes the proof.\endproof
\medskip\noindent
{\it {\bigsll R}}EMARK {\bigsl 3.2.} For $p>2$, Lemma 3.1 was proved by Faltings using Raynaud's theorem [R, 3.3.3] (see [Mo, 3.6] and [V, 3.2.17, Step B]). 
\bigskip\smallskip
\centerline{{\bf \bigsll\bigsll 4. Proof of Theorem 1.3}}
\bigskip 
Let $O$, $K$, $e$, and $Y$ be as in Section 1. We start with a general Proposition. 
\medskip\noindent
{\it {\bigsll P}}ROPOSITION {\bigsl 4.1.} {\it If $Y$ is $p$-healthy regular then $Y$ is also healthy regular.} 
\medskip\noindent
{\it Proof.} Let $(Y,U)$ be an extensible pair with $U$ containing $Y_K$, and let $A_U$ be an abelian scheme over $U$. We need to show that $A_U$ extends to an abelian scheme $A$ over $Y$. Since $Y$ is $p$-healthy regular, the $p$-divisible group $D_U$ of $A_U$ extends to a $p$-divisible group $D$ over $Y$. From now on we forget that $Y$ is $p$-healthy regular and we will use just the existence of $D$ to show that $A$ exists. 
\smallskip
Let $N\in\NN\setminus\{1,2\}$ be prime to $p$. To show that $A$ exists, we can assume that $Y$ is local, complete, and strictly henselian, that $U$ is the complement of the maximal point $y$ of $Y$, and that $A_U$ has a principal polarization $p_{A_U}$ and a level $N$ structure $l_{U,N}$ (see [FaC, (i)-(iii) of pp. 185, 186]). We write $Y={\rm Spec}(R)$. Let $p_{D_U}$ be the principal quasi-polarization of $D_U$ defined naturally by $p_{A_U}$; it extends to a principal quasi-polarization $p_D$ of $D$ (cf. Tate's theorem [T, Thm. 4]). Let $g$ be the relative dimension of $A_U$. Let $\Ma_{g,1,N}$ be the moduli scheme over ${\rm Spec}(\ZZ[{1\over N}])$ parameterizing principally polarized abelian schemes over ${\rm Spec}(\ZZ[{1\over N}])$-schemes, of relative dimension $g$ and with level $N$ structure (see [MFK, 7.9, 7.10]). Let $(\Ma,\Mp_{\Ma})$ be the universal principally polarized abelian scheme over $\Ma_{g,1,N}$. 
\smallskip
Let $f_U:U\to\Ma_{g,1,N}$ be the morphism defined by $(A_U,p_{A_U},l_{U,N})$. We show that $f_U$ extends to a morphism $f_Y:Y\to\Ma_{g,1,N}$. 
\smallskip
 Let $N_0\in\NN$ be prime to $p$. From the classical purity theorem we get that the \'etale cover $A_U[N_0]\to U$ extends to an \'etale cover $Y_{N_0}\to Y$. But as $Y$ is strictly henselian, $Y$ has no connected \'etale cover different from $Y$. So each $Y_{N_0}$ is a disjoint union of $N_0^{2g}$-copies of $Y$. Hence $A_U$ has a level $N_0$ structure $l_{U,N_0}$ for any $N_0\in\NN$ prime to $p$. 
\smallskip
Let $\overline{\Ma}_{g,1,N}$ be a projective, toroidal compactification of $\Ma_{g,1,N}$ such that (a) the complement of $\Ma_{g,1,N}$ in $\overline\Ma_{g,1,N}$ has pure codimension 1 in $\overline\Ma_{g,1,N}$ and (b) there is a semi-abelian scheme over $\overline{\Ma}_{g,1,N}$ extending $\Ma$ (cf. [FaC, Chap. IV, Thm. 6.7]). Let $\tilde Y$ be the normalization of the Zariski closure of $U$ in $Y\times_O {\overline{\Ma}_{g,1,N}}_{O}$. It is a projective, normal, integral $Y$-scheme having $U$ as an open subscheme. Let $C$ be the complement of $U$ in $\tilde Y$ endowed with the reduced structure; it is a reduced, projective scheme over the residue field $k$ of $y$. The $\ZZ$-algebras of global functions of $Y$, $U$ and $\tilde Y$ are all equal to $R$ (cf. [Ma, Thm. 38] for $U$). So $C$ is a connected $k$-scheme (cf. [H, 11.3, p. 279]). 
\smallskip
Let $\overline{A}_{\tilde Y}$ be the semi-abelian scheme over $\tilde Y$ extending $A_{U}$. Owing to existence of the $l_{U,N_0}$'s, the N\'eron--Ogg--Shafarevich criterion (see [BLR, p. 183]) implies that $\overline{A}_{\tilde Y}$ is an abelian scheme in codimension at most 1. Threfore, since the complement of $\Ma_{g,1,N}$ in $\overline\Ma_{g,1,N}$ has pure codimension 1 in $\overline\Ma_{g,1,N}$, it folows that $\overline{A}_{\tilde Y}$ is an abelian scheme. So $f_U$ extends to a morphism $f_{\tilde Y}:\tilde Y\to\Ma_{g,1,N}$. Let $p_{\overline{A}_{\tilde Y}}:=f_{\tilde Y}^*(\Mp_{\Ma})$. Tate's theorem implies that the principally quasi-polarized $p$-divisible group of $(\overline{A}_{\tilde Y},p_{\overline{A}_{\tilde Y}})$ is the pull-back $(D_{\tilde Y},p_{D_{\tilde Y}})$ of $(D,p_D)$ to $\tilde Y$. Hence the pull-back $(D_C,p_{D_C})$ of $(D_{\tilde Y},p_{D\tilde Y})$ to $C$ is constant; that is, it is the pull-back to $C$ of a principally quasi-polarized $p$-divisible group over $k$. 
\smallskip
We check that the image $f_{\tilde Y}(C)$ of $C$ through $f_{\tilde Y}$ is a point $\{y_0\}$ of $\Ma_{g,1,N}$. Since $C$ is connected, to check this it suffices to show that, if $\widehat{O_{c}}$ is the completion of the local ring $O_{c}$ of $C$ at an arbitrary point $c$, then the morphism ${\rm Spec}(\widehat{O_{c}})\to\Ma_{g,1,N}$ defined naturally by $f_{\tilde Y}$ is constant. But as $(D_C,p_{D_C})$ is constant, this follows from Serre--Tate deformation theory  (see [Me, Chaps. 4, 5]). So $f_{\tilde Y}(C)$ is a point $\{y_0\}$ of $\Ma_{g,1,N}$. 
\smallskip
Let $R_0$ be the local ring of $\Ma_{g,1,N}$ at $y_0$. Because $Y$ is local and $\tilde Y$ is a projective $Y$-scheme, each point of $\tilde Y$ specializes to a point of $C$. Hence each point of the image of $f_{\tilde Y}$ specializes to $y_0$ and so $f_{\tilde Y}$ factors through the natural morphism ${\rm Spec}(R_0)\to \Ma_{g,1,N}$. Since $R$ is the ring of global functions of $\tilde Y$, the resulting morphism $\tilde Y\to {\rm Spec}(R_0)$ factors through a morphism ${\rm Spec}(R)\to {\rm Spec}(R_0)$. Therefore, $f_{\tilde Y}$ factors through a morphism $f_Y:Y\to \Ma_{g,1,N}$ extending $f_U$. This ends the argument for the existence of $f_Y$. We conclude that $A:=f_Y^*(\Ma)$ extends $A_U$, which completes the proof.\endproof 
\medskip\noindent
{\it {\bigsll R}}EMARK {\bigsl 4.2.} In the proof of Proposition 4.1, the use of semi-abelian schemes can be replaced by de Jong's good reduction criterion [dJ, 2.5] as follows. If we define $\tilde Y$ to be the normalization of the Zariski closure of $U$ in $Y\times_O {\Ma_{g,1,N}}_O$, then [dJ] implies that the morphism $\tilde Y\to Y$ of $O$-schemes of finite type satisfies the valuative criterion of properness with respect to discrete valuation rings of equal characteristic $p$. Using (as in the proof of Proposition 4.1) the N\'eron--Ogg--Shafarevich criterion, one checks that the morphism $\tilde Y\to Y$ of $O$-schemes satisfies the valuative criterion of properness with respect to discrete valuation rings whose fields of fractions have characteristic $0$. Hence the morphism $\tilde Y\to Y$ of $O$-schemes is proper. The rest of the argument is entirely the same.  
\medskip\noindent
{\it {\bigsll C}}ONCLUSION {\bigsl 4.3.} We assume that $e=1$ and that $Y$ is formally smooth over $O$. Based on Proposition 4.1, in order to prove Theorem 1.3 it suffices to show that $Y$ is $p$-healthy regular. So let $(Y,U)$ be an extensible pair with $U$ containing $Y_K$. We need to show that any $p$-divisible group $D_U$ over $U$ extends to a $p$-divisible group $D$ over $Y$.  This is a local statement for the faithfully flat topology, so we can assume that $Y$ is local, complete, and strictly henselian and that $U$ is the complement of the maximal point $y$ of $Y$ (see [FaC, p. 183]). Write $Y={\rm Spec}(R)$, and let $d\in\NN$ be the dimension of $R/pR$. We show the existence of $D$ by induction on $d$.
\smallskip
If $d=1$ then, for all $n$, $m\in\NN$, the short exact sequence $0\to D_U[p^n]\to D_U[p^{n+m}]\to D_U[p^m]\to 0$ in the category $p-FF(U)$ extends uniquely to a short exact sequence $0\to D_n\to D_{n+m}\to D_m\to 0$ in the category $p-FF(Y)$ (cf. Lemma 3.1). Hence there is a unique $p$-divisible group $D$ over $Y$ such that $D[p^n]=D_n$. Obviously $D$ extends $D_U$. 
For $d\ge 2$, the passage from $d-1$ to $d$ is entirely as in [FaC, pp. 183, 184] applied to $R$ and any regular parameter $x\in R$ such that $R/xR$ is formally smooth over $O$. This ends the induction and so establishes the existence of $D$, concluding the proof of Theorem 1.3. 
\bigskip\smallskip
\centerline{\bigsll {\bf References}}
\bigskip
\item{[Ba]} A. Badra, {\it  D\'eformations des $p$-groupes finis commutatifs sur un corps parfait et filtration de Hodge}, C. R. Acad. Sci. Paris S\'er. A-B 291 (1980), pp. 539--542.
\item{[B]} P. Berthelot, {\it Cohomologie cristalline des sch\'emas de caract\'eristique $p>0$}, Lecture Notes in Math., 407, Springer-Verlag, New York, 1974.
\item{[BBM]} P. Berthelot, L. Breen, and W. Messing, {\it Th\'eorie de Dieudonn\'e cristalline II}, Lecture Notes in Math., 930, Springer-Verlag, New York, 1982.
\item{[BLR]} S. Bosch, W. L\"utkebohmert, and M. Raynaud, {\it N\'eron models}, Springer-Verlag, Berlin, 1990.
\item{[dJ]} J. de Jong, {\it Homomorphisms of Barsotti--Tate groups and crystals in positive characteristic}, Invent. Math. 134 (1998), pp. 301--333.
\item{[dJO]} J. de Jong and F. Oort, {\it On extending families of curves},
J. Algebraic Geom. 6 (1997), pp. 545--562. 
\item{[Fa]} G. Faltings, {\it Crystalline cohomology and $p$-adic Galois representations}, Algebraic analysis, geometry, and number theory (Baltimore, 1998), pp. 25--80, Johns Hopkins Univ. Press, Baltimore, 1989.
\item{[FaC]} G. Faltings and C.-L. Chai, {\it Degeneration of abelian varieties}, Springer-Verlag, Berlin, 1990.
\item{[F]} J.-M. Fontaine, {\it Groupes $p$-divisibles sur les corps locaux}, Ast\'erisque 47/48 (1977).
\item{[FL]} J.-M. Fontaine and G. Laffaille, {\it Construction de repr\'esentations p-adiques}, Ann. Sci. \'Ecole Norm. Sup. (4) 15 (1982), pp. 547--608.
\item{[G]} A. Grothendieck, et al. {\it Rev\^etements \'etales et groupe fondamental}, Lecture Notes in Math., 224, Springer-Verlag, New York, 1971.
\item{[H]} R. Hartshorne, {\it Algebraic geometry}, Grad. Texts in Math., 52, Springer-Verlag, Berlin, 1977.
\item{[Ka]} N. Katz, {\it Travaux de Dwork}, S\'eminaire Bourbaki, 24eme annee (1971/1972), Exp. no. 409, pp. 167-200, Lecture Notes in Math., 317, Springer-Verlag, New York, 1973.
\item{[K]} R. E. Kottwitz, {\it Points on some Shimura varieties over finite fields}, J. Amer. Math. Soc. 5 (1992), pp. 373--444.
\item{[Ma]} H. Matsumura, {\it Commutative algebra}, Benjamin, New York, 1980 (second edition).
\item{[Me]} W. Messing, {\it The crystals associated to Barsotti--Tate groups, with applications to abelian schemes}, Lecture Notes in Math., 264, Springer-Verlag, New York, 1972.
\item{[Mi]} J. S. Milne, {\it The points on a Shimura variety modulo a prime of good reduction}, The Zeta functions of Picard modular surfaces, pp. 153--255, Univ. Montreal Pres, Montreal, Quebec, 1992.
\item{[Mo]} B. Moonen, {\it Models of Shimura varieties in mixed characteristics}, Galois representations in arithmetic
algebraic geometry (Durham, 1996), pp. 267--350, London Math. Soc. Lecture Note Ser., 254, pp. 267--350, Cambridge Univ. Press, Cambridge, U.K., 1998.
\item{[M-B]} L. Moret-Bailly, {\it Un th\'eor\`eme de puret\'e pour les families de courbes lisses}, C. R. Acad. Sci. Paris S\'er. I Math. 300 (1985), pp. 489--492.
\item{[MFK]} D. Mumford, J. Fogarty and F. Kirwan, {\it Geometric invariant theory}, Springer-Verlag, New York, 1994.
\item{[N]} A. N\'eron, {\it Mod\`eles minimaux des vari\'et\'es ab\'eliennes}, Inst. Hautes \'Etudes Sci. Publ. Math. 21 (1964).
\item{[R]} M. Raynaud, {\it Sch\'emas en groupes de type (p,...,p)}, Bull. Soc. Math. France 102 (1974), pp. 241--280.
\item{[T]} J. Tate, {\it p-divisible groups}, Proceedings of a conference on local fields (Driesbergen, 1966), pp. 158--183, Springer-Verlag, Berlin, 1967.
\item{[V]} A. Vasiu, {\it Integral canonical models of Shimura varieties of preabelian
 type}, Asian J. Math. 3 (1999), pp. 401--518.
\vskip 0.5 in
\line{\hfill\vbox{
\vbox{Department of Mathematics}
\vbox{University of Arizona} 
\vbox{617 N. Santa Rita, P.O. Box 210089}
\vbox{Tucson, AZ-85721-0089}
\medskip
\vbox{adrian@math.arizona.edu\ \ \ \ \ }}}

\end